\documentclass[12pt]{article}
\usepackage{cite}
\usepackage{graphicx}
\usepackage[cmex10]{amsmath}
\usepackage[framed,amsmath,thmmarks,hyperref]{ntheorem}
\usepackage{amssymb}
\usepackage{arydshln}
\usepackage{multirow}
\usepackage{caption}
\usepackage[hidelinks]{hyperref}

\newcommand{\R}{\mathbb{R}}

\newtheorem{theorem}{Theorem}
\newtheorem{corollary}{Corollary}
\newtheorem{remark}{Remark}

\begin{document}

\title{Event-Triggered $H_\infty$ Control: a Switching Approach}

\author{Anton~Selivanov \and
        Emilia~Fridman
\thanks{A. Selivanov ({\tt\small antonselivanov@gmail.com}) and E. Fridman ({\tt\small emilia@eng.tau.ac.il}) are with School of Electrical Engineering, Tel Aviv University, Israel.}
\thanks{This work was published in \cite{Selivanov2015f,Selivanov2016c}. Supported by Israel Science Foundation (grant No. 1128/14).}%
}
\renewcommand\footnotemark{}
\maketitle

\begin{abstract}
  Event-triggered approach to networked control systems is used to reduce the workload of the communication network. For the static output-feedback continuous event-trigger may generate an infinite number of sampling instants in finite time (Zeno phenomenon) what makes it inapplicable to the real-world systems. Periodic event-trigger avoids this behavior but does not use all the available information. In the present paper we aim to exploit the advantage of the continuous-time measurements and guarantee a positive lower bound on the inter-event times by introducing a switching approach for finding a waiting time in the event-triggered mechanism. Namely, our idea is to present the closed-loop system as a switching between the system under periodic sampling and the one under continuous event-trigger and take the maximum sampling preserving the stability as the waiting time. We extend this idea to the $L_2$-gain and ISS analysis of perturbed networked control systems with network-induced delays. By examples we demonstrate that the switching approach to event-triggered control can essentially reduce the amount of measurements to be sent through a communication network compared to the existing methods.
\end{abstract}

\section{Introduction}
Networked control systems (NCS), that are comprised of sensors, actuators, and controllers connected through a communication network, have been recently extensively studied by researchers from a variety of disciplines \cite{Antsaklis2004,Hespanha2007,Garcia2013,Fridman2014}. One of the main challenges in such systems is that only sampled in time measurements can be transmitted through a communication network. Namely, consider the system
\begin{equation}\label{LTI}
\dot x(t)=Ax(t)+Bu(t),\quad y(t)=Cx(t),
\end{equation}
with a state $x\in\R^n$, input $u\in\R^m$, and output $y\in\R^l$. Assume that there exists $K\in\R^{m\times l}$ such that the control signal $u(t)=-Ky(t)$ stabilizes the system \eqref{LTI}. In NCS the measurements can be transmitted to the controller only at discrete time instants
\begin{equation}\label{sk}
0=s_0<s_1<s_2<\ldots,\quad \lim_{k\to\infty}s_k=\infty.
\end{equation}
Therefore, the closed-loop system has the form
\begin{equation}\label{CL-NCS}
\dot x(t)=Ax(t)-BKCx(s_k), \quad t\in[s_k,s_{k+1}),\quad k\in\mathbb{N}_0,
\end{equation}
where $\mathbb{N}_0$ is the set of nonnegative integers. There are different ways of obtaining the sequence of sampling instants $s_k$ that preserve the stability. The simplest approach is \textit{periodic sampling} where one chooses $s_k=kh$ with appropriate period~$h$. Under periodic sampling the measurements are sent even when the output fluctuation is small and does not significantly change the control signal. To avoid these ``redundant'' packets one can use \textit{continuous event-trigger} \cite{Tabuada2007}, where
\begin{equation}\label{conEvTr}
s_{k+1}=\min\{t>s_k\,|\,(y(t)-y(s_k))^T\Omega(y(t)-y(s_k))\geq\varepsilon y^T(t)\Omega y(t)\}
\end{equation}
with a matrix $\Omega\geq0$ and a scalar $\varepsilon>0$. In case of a static \textit{output}-feedback execution times $s_k$, implicitly defined by \eqref{conEvTr}, can be such that $\lim_{k\to\infty}s_k<\infty$ \cite{Borgers2014}. That is, an infinite number of events is generated in finite time what makes \eqref{conEvTr} inapplicable to NCS. To avoid this Zeno phenomenon one can use \textit{periodic event-trigger} \cite{Heemels2013a,Peng2013,Yue2013,Zhang2014} by choosing
\begin{equation}\label{tkdir}
s_{k+1}\!=\!\min\{s_k+ih\,|\,i\in\mathbb{N},\,(y(s_k+ih)-y(s_k))^T\Omega
(y(s_k+ih)-y(s_k))\!>\!\varepsilon y^T(s_k+ih)\Omega y(s_k+ih)\}.
\end{equation}
This approach guarantees that the inter-event times are at least~$h$ and fits the case where the sensor measures only sampled in time outputs $y(ih)$.

However, when the continuous measurements are available one can use this additional information to improve the control algorithm. In \cite{Heemels2008,Tallapragada2012,Tallapragada2012b} the following strategy of choosing the sampling instants has been considered:
\begin{equation}\label{T}
s_{k+1}=\min\{t\geq s_k+T\,|\,\eta\geq0\},
\end{equation}
where $T>0$ is a constant waiting time and $\eta$ is an event-trigger condition. In \cite{Tallapragada2012,Tallapragada2012b} the value of $T$ that preserves the stability was obtained by solving a scalar differential equation. For $\eta=|y(t)-y(t_k)|-C$ with a constant $C$ some qualitative results concerning \textit{practical} stability have been obtained in~\cite{Heemels2008}.

In this work we propose a new constructive and efficient method of finding an appropriate waiting time. Our idea is to present the closed-loop system as a switching between the system under periodic sampling and the one under continuous event-trigger and take the maximum sampling preserving the stability as a waiting time. We extend this idea to the systems with network-induced delays, external disturbances, and measurement noise (Section~\ref{Sec:EvTrDel}). Differently from \cite{Heemels2008,Tallapragada2012,Tallapragada2012b,Heemels2013a} our method is applicable to uncertain linear systems and the waiting time is found from LMIs. Comparatively to periodic event-trigger of \cite{Peng2013,Yue2013,Zhang2014} our method leads to error separation between the system under periodic sampling and the one under continuous event-trigger that allows for larger sampling periods for the same values of the event-trigger parameter~$\varepsilon$. The latter allows to reduce the amount of sent measurements as illustrated by examples brought from \cite{Tabuada2007}, \cite{Borgers2014}, and \cite{Wang2009b} (Section~\ref{sec:examples}).

\section{A switching approach to event-trigger}\label{Sec:EvTrNoDel}
Consider \eqref{LTI}. Assume that there exists $K$ such that $A-BKC$ is Hurwitz. For $C=I$ such $K$ exists if $(A,B)$ is stabilizable. For the static output-feedback case such $K$ exists if the transfer function $C(sI-A)^{-1}B$ is hyper-minimum-phase (has stable zeroes and positive leading coefficient of the numerator, see, e.g., \cite{Fradkov2003}). Assume that the measurements are sent at time instants \eqref{sk}. According to \cite{Fridman2004} the closed-loop system \eqref{CL-NCS} under periodic sampling $s_k=kh$ can be presented in the form
\begin{equation}\label{sys1}
\dot x(t)=(A-BKC)x(t)+BKC\int_{t-\tau(t)}^t\dot x(s)\,ds,
\end{equation}
where $\tau(t)=t-s_k$ for $t\in[s_k,s_{k+1})$. The system \eqref{CL-NCS} under continuous event-trigger \eqref{conEvTr} can be rewritten as (see~\cite{Tabuada2007})
\begin{equation}\label{sys2}
\dot x(t)=(A-BKC)x(t)-BKe(t)
\end{equation}
with $e(t)=y(s_k)-y(t)$ for $t\in[s_k,s_{k+1})$.

Under periodic sampling (leading to \eqref{sys1}) ``redundant'' packets can be sent while continuous event-trigger (that leads to \eqref{sys2}) can cause Zeno phenomenon. To avoid the above drawbacks periodic event-trigger \eqref{tkdir} can be used, where the closed-loop system can be written as
\begin{equation}\label{directIdea}
\dot x(t)=(A-BKC)x(t)+BKC\int_{t-\tau(t)}^t\dot x(s)\,ds-BKe(t)
\end{equation}
with $\tau(t)=t-s_k-ih\leq h$, $e(t)=y(s_k)-y(s_k+ih)$ for $t\in[s_k+ih,s_k+(i+1)h)$, $i\in\mathbb{N}_0$ such that $s_k+(i+1)h\leq s_{k+1}$. As one can see, the error due to sampling that appears in \eqref{sys1} (the integral term) and the error $e(t)$ due to triggering from \eqref{sys2} are both presented in \eqref{directIdea} what makes it more difficult to ensure the stability of~\eqref{directIdea} compared to \eqref{sys1} or \eqref{sys2}.

We propose an event-trigger that allows to \textit{separate these errors} by considering the switching between periodic sampling and continuous event-trigger. Namely, after the measurement has been sent, the sensor waits for at least $h$ seconds (that corresponds to $T$ in \eqref{T}). During this time the system is described by \eqref{sys1}. Then the sensor begins to continuously check the event-trigger condition and sends the measurement when it is violated. During this time the system is described by \eqref{sys2}. This leads to the following choice of sampling:
\begin{equation}\label{tk}
s_{k+1}=\min\{s\geq s_k+h\,|\,(y(s)-y(s_k))^T\Omega
(y(s)-y(s_k))\geq\varepsilon y^T(s)\Omega y(s)\}
\end{equation}
with a matrix $\Omega\geq0$ and scalars $\varepsilon\geq0$, $h>0$, where the inter-event times are not less than~$h$. The system \eqref{CL-NCS}, \eqref{tk} can be presented as a switching between \eqref{sys1} and \eqref{sys2}:
\begin{equation}\label{switchSys}
\dot x(t)=(A-BKC)x(t)+\chi(t)BKC\int_{t-\tau(t)}^t\dot x(s)\,ds-(1-\chi(t))BKe(t),
\end{equation}
where
\begin{equation}\label{notations:e-tau}
\begin{aligned}
&\begin{aligned}
\tau(t)&=t-s_k\leq h, &&t\in[s_k,s_k+h),\\
e(t)&=y(s_k)-y(t), &&t\in[s_k+h,s_{k+1}),
\end{aligned}\\
&\chi(t)=\left\{\begin{aligned}
&1,&&t\in[s_k,s_k+h), \\
&0,&&t\in[s_k+h,s_{k+1}).
\end{aligned}\right.
\end{aligned}
\end{equation}

To obtain the stability conditions for the switched system \eqref{switchSys} we use different Lyapunov functions: for \eqref{switchSys} with $\chi(t)=0$ we consider
\begin{equation}\label{VP}
V_P(x)=x^T(t)Px(t), \quad P>0,
\end{equation}
for \eqref{switchSys} with $\chi(t)=1$ we apply the functional from \cite{Fridman2010}:
\begin{equation}\label{LKF}
V(t,x_t,\dot x_t)=V_P(x(t))+V_U(t,\dot x_t)+V_X(t,x_t),
\end{equation}
where $x_t(\theta)=x(t+\theta)$ for $\theta\in[-h,0]$, $V_P$ is given by \eqref{VP},
\begin{equation*}
\begin{aligned}
&V_U(t,\dot x_t)=(h-\tau(t))\int_{t_k}^te^{2\delta(s-t)}\dot x^T(s)U\dot x(s)\,ds,\,U>0, \\
&V_X(t,x_t)=(h-\tau(t))\left[
\begin{matrix}
x(t) \\ x(t_k)
\end{matrix}\right]^T\left[
\begin{matrix}
\frac{X+X^T}2 & -X+X_1 \\
* & -X_1-X_1^T+\frac{X+X^T}2
\end{matrix}\right]
\left[
\begin{matrix}
x(t) \\ x(t_k)
\end{matrix}\right].
\end{aligned}
\end{equation*}
Note that the values of $V$ and $V_P$ coincide at the switching instants $t_k$ and $t_k+h$.

\begin{theorem}\label{th:1}
	For given scalars $h>0$, $\varepsilon\geq0$, $\delta\geq0$ let there exist $n\times n$ matrices $P>0$, $U>0$, $X$, $X_1$, $P_2$, $P_3$, $Y_1$, $Y_2$, $Y_3$ and $l\times l$ matrix $\Omega\geq0$ such that\footnote{MATLAB codes are available at \url{https://github.com/AntonSelivanov/TAC16}}
	\begin{equation}\label{th:1:cond}
	\Xi>0,\quad \Psi_0\leq0,\quad \Psi_1\leq0,\quad\Phi\leq0,
	\end{equation}
	where
	\begin{equation*}
	\begin{aligned}
	\Xi&=\left[\begin{matrix}
	P+h\frac{X+X^T}2 & hX_1-hX \\
	* & -hX_1-hX_1^T+h\frac{X+X^T}2
	\end{matrix}\right],\\
	\Phi&=\left[\begin{matrix}
	\Phi_{11} & \Phi_{12} & -P_2^TBK \\
	* & -P_3^T-P_3 & -P_3^TBK \\
	* & * & -\Omega
	\end{matrix}\right],\\
	\Psi_0&=\left[\begin{matrix}
	\Psi_{11}-X_\delta & \Psi_{12}+h\frac{X+X^T}2 & \Psi_{13}+X_{1\delta} \\
	* & \Psi_{22}+hU & \Psi_{23}-h(X-X_1) \\
	* & * & \Psi_{33}-X_{2\delta}|_{\tau=0}
	\end{matrix}\right],\\
	\Psi_1&=\left[\begin{matrix}
	\Psi_{11}-\frac{X+X^T}2 & \Psi_{12} & \Psi_{13}+X-X_1 & hY_1^T \\
	* & \Psi_{22} & \Psi_{23} & hY_2^T \\
	* & * & \Psi_{33}-X_{2\delta}|_{\tau=h} & hY_3^T \\
	* & * & * & -hUe^{-2\delta h}
	\end{matrix}\right],
	\end{aligned}
	\end{equation*}
	\begin{equation*}
	\begin{array}{l}
	\Phi_{11}=P_2^T(A-BKC)+(A-BKC)^TP_2+\varepsilon C^T\Omega C+2\delta P,\\
	\Phi_{12}=P+(A-BKC)^TP_3-P_2^T,\\
	\Psi_{11}=A^TP_2+P_2^TA+2\delta P-Y_1-Y_1^T, \\
	\Psi_{12}=P-P_2^T+A^TP_3-Y_2,\\
	\Psi_{13}=Y_1^T-P_2^TBKC-Y_3,\\
	\Psi_{22}=-P_3-P_3^T,\\
	\Psi_{23}=Y_2^T-P_3^TBKC,\\
	\Psi_{33}=Y_3+Y_3^T,\quad X_\delta=(1\!/2-\delta h)(X+X^T),\\
	X_{1\delta}=(1-2\delta h)(X-X_1),\\
	X_{2\delta}=(1\!/2-\delta(h-\tau))(X+X^T-2X_1-2X_1^T).
	\end{array}
	\end{equation*}
	Then the system \eqref{CL-NCS} under the event-trigger \eqref{tk} is exponentially stable with a decay rate $\delta$.
\end{theorem}

{\em Proof.}
The system \eqref{CL-NCS}, \eqref{tk} is presented in the form of the switched system \eqref{switchSys}. According to \cite{Fridman2010} the conditions $\Xi>0$, $\Psi_0\leq0$, $\Psi_1\leq0$ imply $V\geq\alpha|x(t)|^2$ and $\dot V\leq-2\delta V$ for the system \eqref{switchSys} with $\chi(t)=1$. Consider \eqref{switchSys} with $\chi(t)=0$. Since for $t\in[s_k+h,s_{k+1})$ the relation \eqref{tk} implies
\begin{equation}\label{EvTr1}
0\leq\varepsilon x^T(t)C^T\Omega Cx(t)-e^T(t)\Omega e(t),
\end{equation}
we add \eqref{EvTr1} to $\dot V_P$ to compensate the cross term with $e(t)$. We have
\begin{multline*}
\dot V_P+2\delta V_P\leq2x^TP\dot x+2\delta x^TPx+2[x^TP_2^T+\dot x^TP_3^T][(A-BKC)x-BKe-\dot x]\\+[\varepsilon x^TC^T\Omega Cx-e^T\Omega e]=\varphi^T\Phi\varphi\leq0,
\end{multline*}
where $\varphi=\operatorname{col}\{x,\dot x,e\}$. Thus, $\dot V_P\leq-2\delta V_P$.

The stability of the switched system \eqref{switchSys} follows from the fact that at the switching instants $s_k$ and $s_k+h$ the values of $V$ and $V_P$ coincide.

\hfill$\square$

By extending the proof from \cite{Fridman2010} we obtain the stability conditions for the system \eqref{CL-NCS}, \eqref{tkdir} presented in the form \eqref{directIdea}:

\begin{remark}\label{pr:1}
For given scalars $h>0$, $\varepsilon\geq0$, $\delta>0$ let there exist $n\times n$ matrices $P>0$, $U>0$, $X$, $X_1$, $P_2$, $P_3$, $Y_1$, $Y_2$, $Y_3$ and $l\times l$ matrix $\Omega\geq0$ such that
\begin{equation}\label{pr:1:cond}
\Xi>0,\quad\left[\begin{array}{c:c}
& -P_2^TBK \\
\overline{\Psi}_i & -P_3^TBK \\
 & 0 \\ \hdashline
 * & -\Omega
\end{array}\right]\leq0,
\end{equation}
where $\overline{\Psi}_i=\Psi_i+\varepsilon[I_n\,0]^TC^T\Omega C[I_n\,0]$, $i=0,1$. Then the system \eqref{CL-NCS} under periodic event-trigger \eqref{tkdir} is exponentially stable with a decay rate $\delta$.
\end{remark}

\begin{remark}\label{rem:1}
The feasibility of \eqref{pr:1:cond} implies the feasibility of~\eqref{th:1:cond}. Therefore, the stability of \eqref{CL-NCS} under \eqref{tk} can be guaranteed for not smaller $h$ and $\varepsilon$ than under~\eqref{tkdir}. Examples in Section~\ref{sec:examples} show that these values under \eqref{tk} are essentially larger what allows to reduce the amount of sent measurements. Note that for \textit{the same} $h$, $\varepsilon$, and $\Omega$ the amount of sent measurements under periodic event-trigger \eqref{tkdir} is deliberately less than under~\eqref{tk}. Indeed, if the measurement is sent at $s_k$ and the event-trigger rule is satisfied at $s_k+h$, according to \eqref{tkdir} the sensor will wait till at least $s_k+2h$ before sending the next measurement, while according to \eqref{tk} the next measurement can be sent before $s_k+2h$.
\end{remark}

\section{Event-trigger under network-induced delays and disturbances}\label{Sec:EvTrDel}
\begin{figure}[!t]
  \centering
  \includegraphics[width=.4\linewidth]{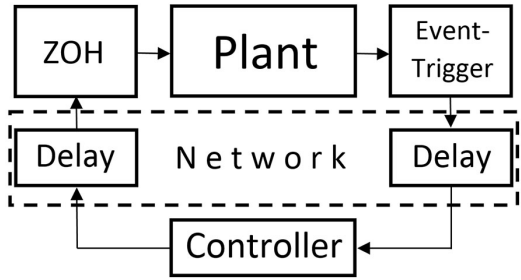}
  \caption{Scheme of a system with network-induced delays}
  \label{fig:scheme-delay}
\end{figure}
Consider the system
\begin{equation}\label{Hinfsys}
\begin{aligned}
\dot x(t)&=Ax(t)+B_1w(t)+B_2u(t), \\
z(t)&=C_1x(t)+D_1u(t), \\
y(t)&=C_2x(t)+D_2v(t)
\end{aligned}
\end{equation}
with a state $x\in\R^n$, input $u\in\R^m$, controlled output $z\in\R^{n_z}$, measurements $y\in\R^l$, and disturbances $w\in\R^{n_w}$, ${v\in\R^{n_v}}$. Denote by $\eta_k\leq\eta_M$ the overall network-induced delay from the sensor to the actuator that affects the transmitted measurement $y(s_k)$ (see Fig.~\ref{fig:scheme-delay}). Here $s_k$ is a sampling instant on the sensor side. We assume that $\eta_k$ are such that the ZOH updating times $t_k=s_k+\eta_k$ satisfy
\begin{equation}\label{etakCond}
t_k=s_k+\eta_k\leq s_{k+1}+\eta_{k+1}=t_{k+1},\quad k\in\mathbb{N}_0.
\end{equation}
Then the system \eqref{Hinfsys} with $u(t)=Ky(s_k)$ for $t\in[t_k,t_{k+1})$ has the form
\begin{equation}\label{sysDelay}
\begin{aligned}
\dot x(t)&=Ax(t)+B_1w(t)+B_2K[C_2x(t_k-\eta_k)+D_2v(t_k-\eta_k)], \\
z(t)&=C_1x(t)+D_1K[C_2x(t_k-\eta_k)+D_2v(t_k-\eta_k)].
\end{aligned}
\end{equation}

Similar to Section~\ref{Sec:EvTrNoDel} we would like to present the resulting closed-loop system \eqref{tk}, \eqref{sysDelay} as a system with periodic sampling for $t\in[t_k,t_k+h)$ (i.e. $t\in[s_k+\eta_k,s_k+\eta_k+h)$) and as a system with continuous event-trigger for $t\in[t_k+h,t_{k+1})$. If $t_k+h=s_k+\eta_k+h>s_{k+1}+\eta_{k+1}=t_{k+1}$ (what may happen due to the communication delay $\eta_k$) no switching occurs. Therefore, the system \eqref{tk}, \eqref{sysDelay} can be presented as
\begin{equation}\label{switchSysDel}
\begin{aligned}
\dot x(t)&=Ax(t)+B_1w(t)+\chi(t)B_2K[C_2x(t-\tau(t))+D_2v(t-\tau(t))]\\
&+(1-\chi(t))B_2K[C_2x(t-\bar{\eta}(t))+D_2v(t-\bar{\eta}(t))+e(t)], \\
z(t)&=C_1x(t)+\chi(t)D_1K[C_2x(t-\tau(t))+D_2v(t-\tau(t))]\\
&+(1-\chi(t))D_1K[C_2x(t-\bar{\eta}(t))+D_2v(t-\bar{\eta}(t))+e(t)],
\end{aligned}
\end{equation}
where
\begin{equation*}
\begin{aligned}
&\chi(t)=\left\{\begin{aligned}
&1,&&t\in[t_k,\min\{t_k+h,t_{k+1}\}), \\
&0,&&t\in[\min\{t_k+h,t_{k+1}\},t_{k+1}),
\end{aligned}\right.\\
&\begin{aligned}
\tau(t)&=t-s_k,&&t\in[t_k,\min\{t_k+h,t_{k+1}\}),\\
e(t)&=y(s_k)-y(t-\bar{\eta}(t)),&&t\in[\min\{t_k+h,t_{k+1}\},t_{k+1}).
\end{aligned}
\end{aligned}
\end{equation*}
Here $\tau(t)\leq h+\eta_M\triangleq\tau_M$ and $\bar{\eta}(t)\in[0,\eta_M]$ is a ``fictitious'' delay to be defined hereafter.

\begin{figure}[!t]
  \centering
  \includegraphics[width=.5\linewidth]{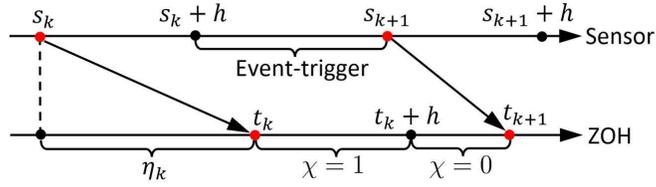}
  \caption{Switching between the subsystems of \eqref{switchSysDel}}
  \label{fig:eta}
\end{figure}

Consider the case where $t_k+h<t_{k+1}$ (see Fig.~\ref{fig:eta}). To use the event-trigger condition we would like to choose such $\bar\eta(t)$ that \eqref{tk} implies
\begin{equation}\label{EvTr}
0\leq\varepsilon[C_2x(t-\bar{\eta}(t))+D_2v(t-\bar{\eta}(t))]^T\Omega[C_2x(t-\bar{\eta}(t))+D_2v(t-\bar{\eta}(t))]-e^T(t)\Omega e(t)
\end{equation}
for $t\in[t_k+h,t_{k+1})$. Relation \eqref{EvTr} is true if $t-\bar\eta(t)\in[s_k+h,s_{k+1})$ for $t\in[t_k+h,t_{k+1})$. Therefore, the simplest choice of $\bar\eta(t)$ is a linear function with $\bar\eta(t_k+h)=\eta_k$ and $\bar\eta(t_{k+1})=\eta_{k+1}$, i.e. for $t\in[\min\{t_k+h,t_{k+1}\},t_{k+1})$
$$
\bar\eta(t)=\frac{t_{k+1}-t}{t_{k+1}-t_k-h}\eta_k+\frac{t-t_k-h}{t_{k+1}-t_k-h}\eta_{k+1}.
$$

Though for both $\chi(t)=0$ and $\chi(t)=1$ the system \eqref{switchSysDel} includes time-delays, the upper bound $\eta_M$ for $\bar\eta(t)$ is smaller than $\tau_M$ since $\tau(t)$ includes the delay due to sampling.

Define $\tau(t)=\bar{\eta}(t)$ for $t\in[\min\{t_k+h,t_{k+1}\},t_{k+1})$. We say that the system \eqref{tk}, \eqref{sysDelay} has an $L_2$-gain ($H_\infty$ gain) less than $\gamma$ if for the zero initial condition $x(0)=0$ and all $w,v\in L_2[0,\infty)$ such that $w^T(t)w(t)+v^T(t-\tau(t))v(t-\tau(t))\not\equiv0$ the following relation holds on the trajectories of~\eqref{tk},~\eqref{sysDelay}:
\begin{equation}\label{L2gain}
J=\int_0^\infty\Bigl\{z^T(t)z(t)-\gamma^2[w^T(t)w(t)+v^T(t-\tau(t))v(t-\tau(t))]\Bigr\}\,dt<0.
\end{equation}

\begin{theorem}\label{th:2}
For given $\gamma>0$, $h>0$, $\eta_M\geq0$, $\varepsilon\geq0$, $\delta>0$ let there exist $n\times n$ matrices $P>0$, $S_0\geq0$, $S_1\geq0$, $R_0\geq0$, $R_1\geq0$, $G_1$, $G_0$ and $l\times l$ matrix $\Omega\geq0$ such that
\begin{equation}\label{th:2:cond}
\Psi\leq0, \quad \Phi\leq0, \quad \left[\begin{matrix}
R_0 & G_0\\
G_0^T & R_0
\end{matrix}\right]\geq0, \quad
\left[\begin{matrix}
R_1 & G_1\\
G_1^T & R_1
\end{matrix}\right]\geq0,
\end{equation}
where $\Psi=\{\Psi_{ij}\}$ and $\Phi=\{\Phi_{ij}\}$ are symmetric matrices composed from the matrices
\begin{equation*}
\begin{array}{l}
\Psi_{11}=\Phi_{11}=A^TP+PA+2\delta P+S_0-e^{-2\delta\eta_M}R_0+C_1^TC_1,\\
\Psi_{12}=e^{-2\delta\eta_M}R_0,\\
\Psi_{14}=PB_2KC_2+C_1^TD_1KC_2,\\
\Psi_{15}=\Phi_{16}=PB_1,\\
\Psi_{16}=\Phi_{17}=PB_2KD_2+C_1^TD_1KD_2,\\
\Psi_{17}=\Phi_{18}=A^TH,\\ \Psi_{22}=\Phi_{22}=e^{-2\delta\eta_M}(S_1-S_0-R_0)-e^{-2\delta\tau_M}R_1,\\
\Psi_{23}=e^{-2\delta\tau_M}G_1,\\
\Psi_{24}=e^{-2\delta\tau_M}(R_1-G_1),\\
\Psi_{33}=\Phi_{33}=-e^{-2\delta\tau_M}(R_1+S_1),\\
\Psi_{34}=e^{-2\delta\tau_M}(R_1-G_1^T), \\
\Psi_{44}=e^{-2\delta\tau_M}(G_1+G_1^T-2R_1)+(D_1KC_2)^TD_1KC_2,\\
\Psi_{46}=(D_1KC_2)^TD_1KD_2,\\
\Psi_{47}=\Phi_{48}=(B_2KC_2)^TH,\,\Psi_{55}=\Phi_{66}=-\gamma^2I,\\
\Psi_{57}=\Phi_{68}=B_1^TH,\\
\Psi_{77}=\Phi_{88}=-H,\\
\end{array}
\end{equation*}
\begin{equation*}
\begin{array}{l}
\Psi_{66}=(D_1KD_2)^TD_1KD_2-\gamma^2I,\\
\Psi_{67}=\Phi_{78}=(B_2KD_2)^TH,\\
\Phi_{12}=e^{-2\delta\eta_M}G_0,\\
\Phi_{14}=PB_2KC_2+e^{-2\delta\eta_M}(R_0-G_0)+C_1^TD_1KC_2,\\
\Phi_{23}=e^{-2\delta\tau_M}R_1,\\
\Phi_{24}=e^{-2\delta\eta_M}(R_0-G_0^T),\\
\Phi_{15}=PB_2K+C_1^TD_1K,\\
\Phi_{44}=e^{-2\delta\eta_M}(G_0+G_0^T-2R_0)+(D_1KC_2)^TD_1KC_2+\varepsilon C_2^T\Omega C_2,\\
\Phi_{45}=(D_1KC_2)^TD_1K,\\
\Phi_{47}=(D_1KC_2)^TD_1KD_2+\varepsilon C_2^T\Omega D_2,\\
\Phi_{55}=(D_1K)^TD_1K-\Omega,\\
\Phi_{57}=(D_1K)^TD_1KD_2,\\
\Phi_{58}=(B_2K)^TH,\\
\Phi_{77}=(D_1KD_2)^TD_1KD_2+\varepsilon D_2^T\Omega D_2-\gamma^2I,\\
H=\eta_M^2R_0+h^2R_1,
\end{array}
\end{equation*}
$\tau_M=h+\eta_M$, other blocks are zero matrices. Then the system \eqref{sysDelay} under the event-trigger \eqref{tk} is internally exponentially stable with a decay rate $\delta$ and has $L_2$-gain less than $\gamma$.
\end{theorem}

\noindent{\em Proof:} See Appendix.

\begin{corollary}
If \eqref{th:2:cond} are valid with $C_1=0$, $D_1=0$ then the system \eqref{switchSysDel} under the event-trigger \eqref{tk} is Input-to-State Stable with respect to $\bar w(t)=\operatorname{col}\{w(t),v(t-\tau(t))\}$.
\end{corollary}
{\em Proof.}
If $\bar{w}^T(t)\bar{w}(t)$ is bounded by $\Delta^2$ then \eqref{dotV} (see Appendix) with $C_1=0$, $D_1=0$ transforms to $\dot{V}\leq-2\delta V+\gamma^2\Delta^2$. This implies the assertion of the corollary.

\hfill$\square$

\begin{remark}\label{delPar}
The system \eqref{sysDelay} under periodic event-trigger \eqref{tkdir} can be presented in the form \eqref{switchSysDel} with $\chi=0$ and $\bar{\eta}(t)\leq\tau_M$. By modifying the proof of Theorem~\ref{th:2} one can obtain the stability conditions using the functional \eqref{fun} with arbitrary chosen ``delay partitioning'' parameter $\eta_M\in(0,\tau_M)$ \cite{Fridman2009b,Guan2003}.
\end{remark}

\begin{remark}
The proposed approach can take into account packet dropouts with bounded amount of consecutive packet losses and acknowledgement signal of successful reception as suggested in, e.g., \cite{Guinaldo2012}.
\end{remark}


\begin{remark}
Differently from periodic event-trigger approach considered in \cite{Heemels2013a} our method is applicable to linear systems with polytopic-type uncertainties, since LMIs of Theorems~\ref{th:1} and \ref{th:2} are affine in~$A$, $B$, $B_1$, and $B_2$.
\end{remark}

\begin{remark}
	MATLAB codes for solving the LMIs of Theorems~\ref{th:1}, \ref{th:2}, Remarks~\ref{pr:1}, \ref{delPar} are available at \url{https://github.com/AntonSelivanov/TAC16}. 
\end{remark}

\section{Numerical examples}\label{sec:examples}
\subsection*{Example 1 \cite{Tabuada2007}}
Consider the system \eqref{CL-NCS} with
\begin{equation*}
A=\begin{bmatrix}
0 & 1 \\ -2 & 3
\end{bmatrix},\quad
B=\begin{bmatrix}
0 \\ 1
\end{bmatrix},\quad
C=I, \quad
K=\begin{bmatrix}
-1 & 4
\end{bmatrix}.
\end{equation*}
For $\varepsilon=0$ \eqref{tk} transforms into periodic sampling, therefore, Theorem~\ref{th:1} can be used to obtain the maximum period $h$. Under periodic sampling the amount of sent measurements is $\left[\frac{T_f}{h}\right]+1$, where $T_f$ is the time of simulation and $[\cdot]$ is the integer part of a given number. To obtain the amount of sent measurements for $t_k$ given by \eqref{tkdir} (or \eqref{tk}), for each $\varepsilon=i\times10^{-4}$ ($i=0,\ldots,10^{4}$) we find the maximum $h$ that satisfies the conditions of Remark~\ref{pr:1} (or Theorem~\ref{th:1}) and for each pair of $(\varepsilon,h)$ we perform numerical simulations with $T_f=20$ for several initial conditions given by
\begin{equation*}
(x_1(0),x_2(0))=\left(10\cos\left(\frac{2\pi}{30}\,k\right),10\sin\left(\frac{2\pi}{30}\,k\right)\right)
\end{equation*}
with $k=1,\ldots,30$. Then we choose the pair $(\varepsilon,h)$ that ensures the minimum average amount of sent measurements. In this example the best result was achieved under periodic sampling ($\varepsilon=0$). Theorem~\ref{th:1} gives $h=0.356$ for $\delta=0.24$ and $h=0.424$ for $\delta=0.001$. Both event-triggers \eqref{tkdir} and \eqref{tk} did not succeed in reducing the network workload.

\subsection*{Example 2 \cite{Borgers2014}}
Consider the system \eqref{CL-NCS} with
\begin{equation}\label{sysPar}
A=\left[\begin{matrix}
0 & 1 \\ 0 & -3
\end{matrix}\right],\quad
B=\left[\begin{matrix}
0 \\ 1
\end{matrix}\right],\quad
C=\left[\begin{matrix}
1 & 0
\end{matrix}\right], \quad
K=3.
\end{equation}

\begin{table}[!t]
  \renewcommand{\arraystretch}{1}
  \caption{Example 2. Average amounts of sent measurements (SM)}
  \label{Table:noDelay}
  \centering
  \begin{tabular}{lccc}
  & $\varepsilon$ & $h$ & SM \\ \hline
  Periodic sampling & --- & $1.173$ & $18$ \\
  Event-trigger \eqref{tkdir} & $4.6\times10^{-3}$ & $1.115$ & $17.47$ \\
  Event-trigger \eqref{tkdir} & $0.555$ & $0.344$ & $24.8$ \\
  Switching approach \eqref{tk} & $0.555$ & $0.899$ & $11.13$ \\
  \end{tabular}
\end{table}

\begin{table}
  \renewcommand{\arraystretch}{1}
  \caption{Example 2. Average amounts of sent measurements (SM) for different $\eta_M$}
  \label{Table:Delay}
  \centering
  \begin{tabular}{lcccccccc}
  & $\eta_M$ & $0.1$ & $0.2$ & $0.4$ & $0.6$ & $0.7$ \\ \hline
  \multirow{2}{*}{\hspace{-.3cm}\begin{tabular}{l}
  Period. samp./\\
  event-tr. \eqref{tkdir}
  \end{tabular}} & h & $0.636$ & $0.548$ & $0.355$ & $0.143$ & $0.025$ \\
  & SM & $33$  & $38$ & $57.33$ & $139.27$ & $785.73$ \\ \hline
  \multirow{3}{*}{\hspace{-.3cm}\begin{tabular}{l} Event- \\trigger \eqref{tk}\end{tabular}} & $\varepsilon$ & $0.56$ & $0.345$ & $0.075$ & $0.005$ & $0$ \\
  & $h$ & $0.339$ & $0.379$ & $0.278$ & $0.12$ & $0.025$ \\
  & SM & $23.7$ & $28.5$ & $52.4$ & $137.77$ & $785.73$ \\
  \end{tabular}
\end{table}

As it has been shown in \cite{Borgers2014} for this system an accumulation of events occurs under continuous event-trigger \eqref{conEvTr}. In what follows we compare three approaches of choosing the sampling instants~$s_k$: periodic sampling with $s_k=kh$, periodic event-trigger \eqref{tkdir}, and switching event-trigger \eqref{tk}.

\begin{figure}[!t]
  \centering
  \begin{minipage}{.45\textwidth}
    \centering
    \includegraphics[width=\linewidth]{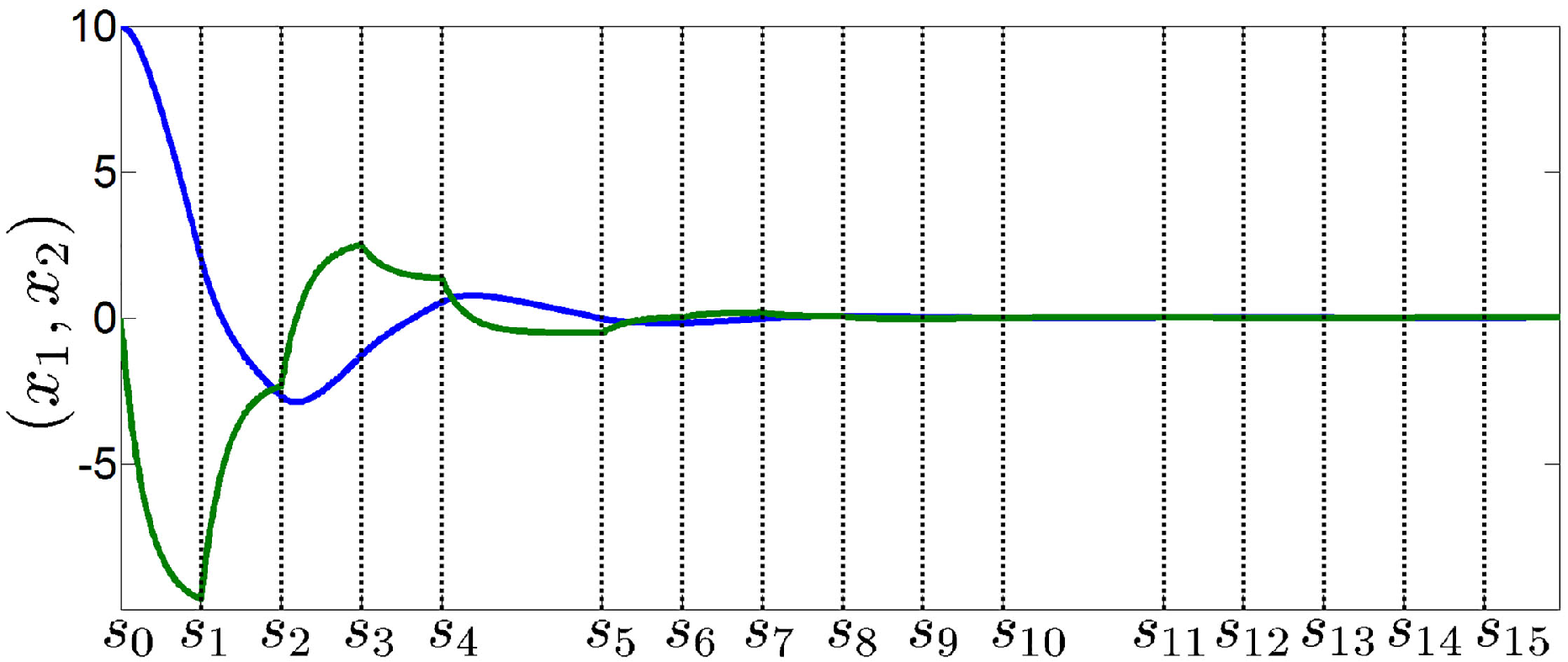}
    \caption{Example 2. Event-trigger \eqref{tkdir}: simulation of the system \eqref{CL-NCS},~\eqref{sysPar}, where $\varepsilon=4.6\times10^{-3}$, $h=1.115$, $[x_1(0),x_2(0)]=[10, 0]$ ($\eta_M=0$).}
    \label{fig:sampling_and_EvTr}
  \end{minipage}%
  \qquad
  \begin{minipage}{.45\textwidth}
    \centering
    \includegraphics[width=\linewidth]{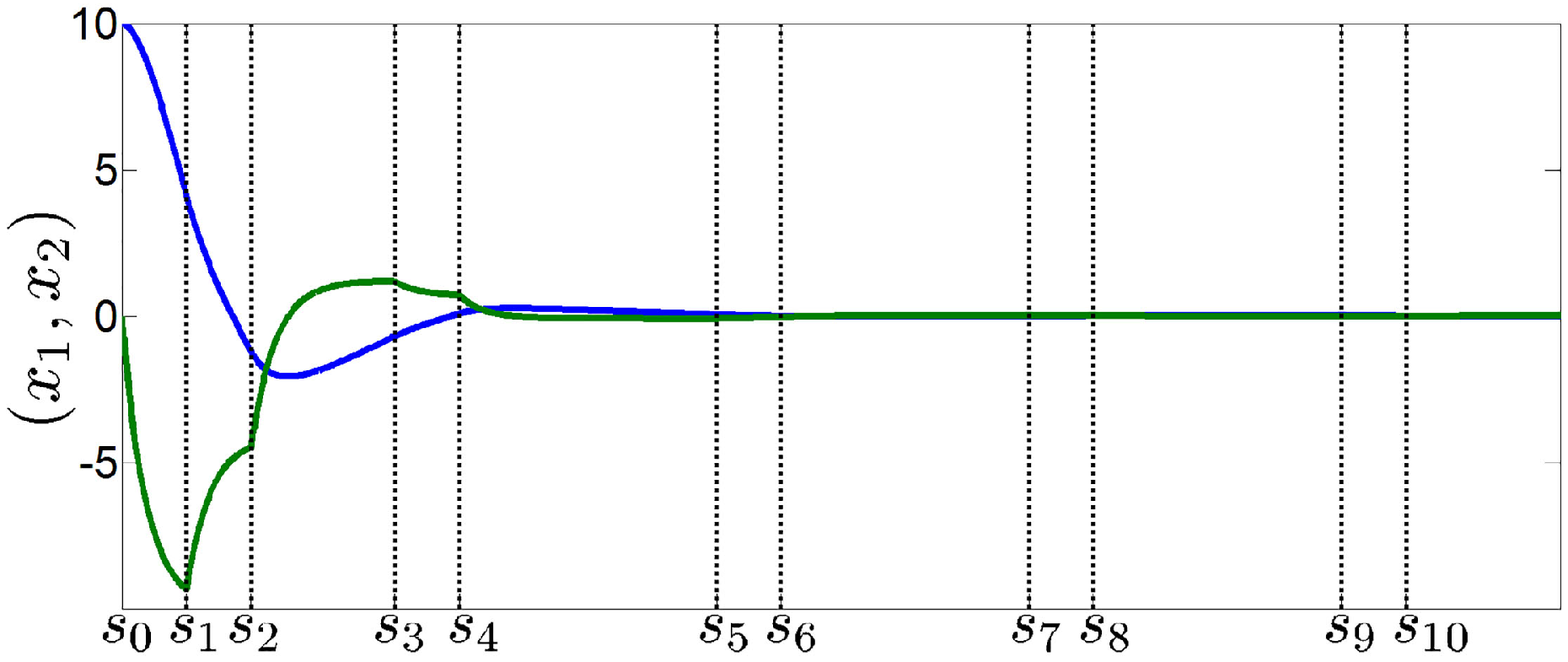}
    \caption{Example 2. Event-trigger \eqref{tk}: simulation of the system \eqref{CL-NCS},~\eqref{sysPar}, where $\varepsilon=0.555$, $h=0.899$, $[x_1(0),x_2(0)]=[10, 0]$ ($\eta_M=0$).}
    \label{fig:Switching}
  \end{minipage}
\end{figure}

We obtained the amount of sent measurements as described in Example~1 (taking $\delta=0.24$, $T_f=20$). As one can see from Table~\ref{Table:noDelay} periodic event-triggered \eqref{tkdir} does not give any significant improvement compared to periodic sampling, while the event-trigger \eqref{tk} allows to reduce the average amount of sent measurements by almost $40\%$. In Figs.~\ref{fig:sampling_and_EvTr}~ and~\ref{fig:Switching} one can see the results of numerical simulations for the event-triggers \eqref{tkdir} and \eqref{tk}. The vertical lines correspond to the time instants when the measurements are sent. The event-trigger \eqref{tkdir} allows to skip the sending of two measurements (after $t_4$ and $t_{10}$), while \eqref{tk} results in large inter sampling times $[t_2,t_3]$, $[t_4,t_5]$,~etc. This allows to significantly reduce the network workload while the decay rate of convergence is preserved.

Now we study the system \eqref{CL-NCS} under network delays. According to the numerical simulations periodic event-trigger \eqref{tkdir} does not give any improvement compared to the periodic sampling for any choice of $\eta_M$ (Remark~\ref{delPar}). Using Theorem~\ref{th:2} with $B_2=B$, $C_2=C$ and other matrices equal to zero we obtained the values of $h$ and $\varepsilon$ in a manner similar to the previously described one. The delays $\eta_k\leq\eta_M$ have been chosen randomly subject to \eqref{etakCond}. The values of $\varepsilon$ and $h$ for $\delta=0.24$ and the corresponding average amounts of sent measurements (SM) during $20$ seconds of simulations for different maximum allowable delays $\eta_M$ are given in Table~\ref{Table:Delay}. As one can see the reduction in the amount of sent measurements vanishes when $\eta_M$ gets larger. This is due to the fact that with the increase of $\eta_M$ the sampling $h$ that preserves the stability is getting smaller, therefore, the difference between $\tau(t)$ and $\bar\eta(t)$ in \eqref{switchSysDel} is getting less significant and the error separation principle proposed here loses its efficiency. However, for $\eta_M=0.1$ the switching approach \eqref{tk} reduces the average amount of the sent measurements by almost $20\%$ while the decay rate of convergence $\delta$ is preserved.

\subsection*{Example 3 \cite{Wang2009b}}
\begin{table}
  \caption{Example 3. Amounts of sent measurements (SM) with time-delays and disturbances}
  \label{Table:delay}
  \centering
  \begin{tabular}{lccc}
  & $\varepsilon$ & $h$ & SM \\ \hline
  Periodic sampling & --- & $0.091$ & $330$ \\
  Event-trigger \eqref{tkdir} & $0.033$ & $0.036$ & $195$ \\
  Event-trigger \eqref{tk} & $0.044$ & $0.065$ & $173$ \\
  \end{tabular}
\end{table}
Consider an inverted pendulum on a cart described by \eqref{CL-NCS} with
\begin{equation*}
A=\left[\begin{matrix}
0 & 1 & 0 & 0 \\
0 & 0 & -1 & 0 \\
0 & 0 & 0 & 1 \\
0 & 0 & 10/3 & 0
\end{matrix}\right],\quad
B=\left[\begin{matrix}
0 \\ 0.1 \\ 0 \\ -1/30
\end{matrix}\right],\quad C=I.
\end{equation*}
For $K=-[2, 12, 378, 210]$ Theorem~\ref{th:1} gives $h=0.242$, $\varepsilon=0.35$. According to the numerical simulations, performed for $T_f$ and $x(0)$ from \cite{Wang2009b}, the average release period under switching event-trigger \eqref{tk} is $0.5769$, which is larger than $0.5131$ obtained for the same system in \cite{Yue2013} (where the average release period is larger than in \cite{Wang2009c,Wang2009b,Carnevale2007,Tabuada2006}).

Consider the system \eqref{sysDelay} with the same $A$, $B_1^T=C_1=[1, 1, 1, 1]$, $B_2=B$, $C_2=I$, $D_1=0.1$, $D_2=[0, 0, 0, 0]^T$, $K=[2.9129, 10.4357, 287.9029, 160.3271]$. For $\gamma=200$, $\eta_M=0.1$ Theorem~\ref{th:2} (with $\delta=0$) gives $h=0.117$, $\varepsilon=0.13$. From the numerical simulations, performed for $T_f$ and $w(t)$ from \cite{Yue2013}, we obtained an average release period $0.3488$, which is larger than $0.3098$ obtained for the same system in \cite{Yue2013} (where the average release period is larger than the one obtained in \cite{Wang2009b} for a different controller gain).

For $\gamma=100$ in a manner similar to Example~1 we obtained the amount of sent measurements presented in Table~\ref{Table:delay}. As one can see both event-triggers reduce the network workload and switching event-trigger \eqref{tk} allows to reduce the amount of sent measurements by more than $11\%$ compared to periodic event-trigger \eqref{tkdir}.

\section{Conclusion}
We proposed a new approach to event-triggered control under the continuous-time measurements that guarantees a positive lower bound for inter-event times and can significantly reduce the workload of the network. Our idea is based on a switching between periodic sampling and continuous event-trigger. We extended this approach to the $L_2$-gain and ISS analyses of perturbed NCS with network-induced delays. Our results are applicable to linear systems with polytopic-type uncertainties. The presented method can be extended to nonlinear NCSs that may be a topic for the future research.

\section*{Appendix\\Proof of Theorem~\ref{th:2}}
The system \eqref{tk}, \eqref{sysDelay} is rewritten as \eqref{switchSysDel}. Similar to \cite{Fridman2009b} we consider Lyapunov functional
\begin{equation}\label{fun}
V=V_P+V_{S_0}+V_{S_1}+V_{R_0}+V_{R_1},
\end{equation}
where $x_t(\theta)=x(t+\theta)$ for $\theta\in[-h,0]$, $V_P(x_t)=x^T(t)Px(t)$,
\begin{equation*}
\begin{aligned}
V_{S_0}(t,x_t)&=\int_{t-\eta_M}^te^{2\delta(s-t)}x^T(s)S_0x(s)\,ds,\\
\end{aligned}
\end{equation*}
\begin{equation*}
\begin{aligned}
V_{R_0}(t,x_t)&=\eta_M\int_{-\eta_M}^0\int_{t+\theta}^te^{2\delta(s-t)}\dot x^T(s)R_0\dot x(s)\,ds\,d\theta, \\
V_{S_1}(t,x_t)&=\int_{t-\tau_M}^{t-\eta_M}e^{2\delta(s-t)}x^T(s)S_1x(s)\,ds,\\
V_{R_1}(t,x_t)&=h\int_{-\tau_M}^{-\eta_M}\int_{t+\theta}^te^{2\delta(s-t)}\dot x^T(s)R_1\dot x(s)\,ds\,d\theta.
\end{aligned}
\end{equation*}
By differentiating these functionals we obtain
\begin{equation}\label{dotVCom}
\begin{aligned}
\dot V_{S_0}&=-2\delta V_{S_0}+x^T(t)S_0x(t)-e^{-2\delta\eta_M}x^T(t-\eta_M)S_0x(t-\eta_M),\\
\dot V_{S_1}&=-2\delta V_{S_1}+e^{-2\delta\eta_M}x^T(t-\eta_M)S_1x(t-\eta_M)-e^{-2\delta\tau_M}x^T(t-\tau_M)S_1x(t-\tau_M), \\
\dot V_{R_0}&=-2\delta V_{R_0}+\eta_M^2\dot x^T(t)R_0\dot x(t)-\eta_M\int_{t-\eta_M}^te^{2\delta(s-t)}\dot x^T(s)R_0\dot x(s)\,ds,\\
\dot V_{R_1}&=-2\delta V_{R_1}+h^2\dot x^T(t)R_1\dot x(t)-h\int_{t-\tau_M}^{t-\eta_M}e^{2\delta(s-t)}\dot x^T(s)R_1\dot x(s)\,ds.
\end{aligned}
\end{equation}
\textit{A. System \eqref{switchSysDel} with $\chi(t)=0$, $\bar\eta(t)\in[0,\eta_M]$.} We have
\begin{equation}\label{dotVP2}
\dot V_P=2x^T(t)P[Ax(t)+B_1w(t)+B_2KC_2x(t-\bar{\eta}(t))+B_2KD_2v(t-\bar\eta(t))+B_2Ke(t)].
\end{equation}
To compensate $x(t-\bar{\eta}(t))$ we apply Jensen's inequality \cite{Gu2003} and Park's theorem \cite{Park2011} to obtain
\small
\begin{equation}\label{eq3}
-\eta_M\int_{t-\eta_M}^te^{2\delta(s-t)}\dot x^T(s)R_0\dot x(s)\,ds
\leq-e^{-2\delta\eta_M}\left[\begin{smallmatrix}
x(t)-x(t-\bar{\eta}(t))\\
x(t-\bar{\eta}(t))-x(t-\eta_M)
\end{smallmatrix}\right]^T
\left[\begin{smallmatrix}
R_0 & G_0\\
G_0^T & R_0
\end{smallmatrix}\right]
\left[\begin{smallmatrix}
x(t)-x(t-\bar{\eta}(t))\\
x(t-\bar{\eta}(t))-x(t-\eta_M)
\end{smallmatrix}\right],
\end{equation}
\begin{equation}\label{eq4}
-h\int_{t-\tau_M}^{t-\eta_M}e^{2\delta(s-t)}\dot x^T(s)R_1\dot x(s)ds
\leq-e^{-2\delta\tau_M}[x(t-\eta_M)-x(t-\tau_M)]^TR_1[x(t-\eta_M)-x(t-\tau_M)].
\end{equation}
\normalsize
By summing up \eqref{EvTr}, \eqref{dotVCom}, \eqref{dotVP2} in view of \eqref{eq3} and \eqref{eq4} and substituting $z$ from \eqref{switchSysDel} we obtain
\begin{equation*}
\dot V+2\delta V+z^Tz-\gamma^2[w^Tw+v^T(t-\bar\eta(t))v(t-\bar\eta(t))]\leq\varphi^T(t)\Phi'\varphi(t)+\dot x^T(t)H\dot x(t),
\end{equation*}
where $\varphi(t)=\operatorname{col}\{x(t), x(t-\eta_M), x(t-\tau_M), x(t-\bar{\eta}(t)),e(t), w(t),v(t-\bar\eta(t))\}$ and the matrix $\Phi'$ is obtained from $\Phi$ by deleting the last block-column and the last block-row. Substituting expression for $\dot x$ and applying Schur complement formula we find that $\Phi\leq0$ guarantees that
\begin{equation}\label{dotV}
\dot V+2\delta V+z^Tz-\gamma^2[w^Tw+v^T(t-\tau(t))v(t-\tau(t))]\leq0.
\end{equation}
\textit{B. System \eqref{switchSysDel} with $\chi=1$, $\tau(t)\in(\eta_M,\tau_M]$.} For $\tau(t)\in[0,\eta_M]$ the system \eqref{switchSysDel} with $\chi=1$ is described by \eqref{switchSysDel} with $\chi=0$ and $e(t)=0$ satisfying \eqref{EvTr}. That is, $\Phi\leq0$ guarantees \eqref{dotV} for \eqref{switchSysDel} with $\chi=1$, $\tau(t)\in[0,\eta_M]$. Therefore, we study the system \eqref{switchSysDel} for $\chi=1$, $\tau(t)\in(\eta_M,\tau_M]$. We have
\begin{equation}\label{dotVP1}
\dot V_P=2x^T(t)P[Ax(t)+B_1w(t)+B_2KC_2x(t-\tau(t))+B_2KD_2v(t-\tau(t))].
\end{equation}
To compensate $x(t-\tau(t))$ with $\tau(t)\in(\eta_M,\tau_M]$ we apply Jensen's inequality and Park's theorem to obtain
\small
\begin{equation}\label{eq1}
-\eta_M\int_{t-\eta_M}^te^{2\delta(s-t)}\dot x^T(s)R_0\dot x(s)\,ds
\leq-e^{-2\delta\eta_M}[x(t)-x(t-\eta_M)]^TR_0[x(t)-x(t-\eta_M)],
\end{equation}
\begin{equation}\label{eq2}
-h\int_{t-\tau_M}^{t-\eta_M}e^{2\delta(s-t)}\dot x^T(s)R_1\dot x(s)ds
\leq-e^{-2\delta\tau_M}\left[\begin{smallmatrix}
x(t-\eta_M)-x(t-\tau(t))\\
x(t-\tau(t))-x(t-\tau_M)
\end{smallmatrix}\right]^T
\left[\begin{smallmatrix}
R_1 & G_1\\
G_1^T & R_1
\end{smallmatrix}\right]
\left[\begin{smallmatrix}
x(t-\eta_M)-x(t-\tau(t))\\
x(t-\tau(t))-x(t-\tau_M)
\end{smallmatrix}\right].
\end{equation}
\normalsize
By summing up \eqref{dotVCom} and \eqref{dotVP1} in view of \eqref{eq1} and \eqref{eq2} and substituting $z$ from \eqref{switchSysDel} we obtain
\begin{equation*}
\dot V+2\delta V+z^Tz-\gamma^2[w^Tw+v^T(t-\tau(t))v(t-\tau(t))]\leq\psi^T(t)\Psi'\psi(t)+\dot x^T(t)H\dot x(t),
\end{equation*}
where $\psi(t)=\operatorname{col}\{x(t), x(t-\eta_M), x(t-\tau_M), x(t-\tau(t)), w(t), v(t-\tau(t))\}$ and the matrix $\Psi'$  is obtained from $\Psi$ by deleting the last block-column and the last block-row. Substituting expression for $\dot x$ and applying Schur complement formula we find that $\Psi\leq0$ guarantees \eqref{dotV} for~\eqref{switchSysDel} with $\chi=1$.

Thus, \eqref{dotV} is true for the switched system \eqref{switchSysDel}. For $w\equiv0$, $v\equiv0$ \eqref{dotV} implies $\dot V\leq-2\delta V$. Therefore, the system \eqref{switchSysDel} is internally exponentially stable with the decay rate $\delta$. By integrating \eqref{dotV} from $0$ to $\infty$ with $x(0)=0$ we obtain~\eqref{L2gain}.

\bibliographystyle{IEEEtran}
\bibliography{IEEEabrv,library}

\end{document}